\documentclass[11pt]{article}
\usepackage{graphicx}
\usepackage{amsmath}
\usepackage{amsthm}
\usepackage{amssymb}
\usepackage{amsfonts}
\usepackage{amsxtra}
\usepackage{epsfig}
\usepackage{verbatim}
\usepackage{bm}
\usepackage{subfigure}
\usepackage{algorithm}
\usepackage{cleveref}
\usepackage[noend]{algpseudocode}
\usepackage{authblk}
\usepackage{algorithmicx}

\makeatletter
\def\BState{\State\hskip-\ALG@thistlm}
\makeatother

\begin{document}

\title{A kernel-independent treecode 
based on barycentric Lagrange interpolation}
\author[1]{L. Wang}
\author[2]{R. Krasny}
\author[3]{S. Tlupova}
\affil[1]{University of Wisconsin-Milwaukee (wang256@uwm.edu)}
\affil[2]{University of Michigan-Ann Arbor (krasny@umich.edu)}
\affil[3]{Farmingdale State College (tlupovs@farmingdale.edu)}
\date{October 22, 2019}
\maketitle

\begin{abstract}
A kernel-independent treecode (KITC) is presented for fast summation of 
particle interactions.
The method employs barycentric Lagrange interpolation at Chebyshev points
to approximate well-separated particle-cluster interactions.
The KITC requires only kernel evaluations,
is suitable for non-oscillatory kernels,
and
it utilizes a scale-invariance property of barycentric Lagrange interpolation.
For a given level of accuracy,
the treecode reduces the operation count for 
pairwise interactions
from $O(N^2)$ to $O(N\log N)$,
where $N$ is the number of particles in the system.
The algorithm is demonstrated 
for systems of regularized Stokeslets and rotlets in 3D,
and
numerical results show the treecode performance 
in terms of error, CPU time, and memory consumption.
The KITC is a relatively simple algorithm with low memory consumption,
and
this enables a straightforward OpenMP parallelization.
\end{abstract}



\section{Introduction}

Consider the problem of evaluating the sum
\begin{equation}
\label{N-body}
u({\bf x}_i) = \sum_{j=1}^N k({\bf x}_i,{\bf x}_j) f_j, \quad i=1,\ldots,N,
\end{equation}
where $u({\bf x}_i)$ is a velocity (or potential or force)
and
$\{{\bf x}_i\} \subset \mathbb{R}^d$ is a set of particles with weights $\{f_i\}$.
Depending on the application, the velocity and weights may be scalars or vectors,
and the kernel may be a tensor.
The kernel $k({\bf x},{\bf y})$ describes the interaction between
a target particle ${\bf x}$ and a source particle ${\bf y}$,
and
we are interested in non-oscillatory kernels that are smooth for ${\bf x} \ne {\bf y}$
and
decay slowly for $|{\bf x}-{\bf y}| \to \infty$.
It is understood that if the kernel is singular for ${\bf x} = {\bf y}$,
then the sum omits the $i=j$ term.

These types of sums arise in particle simulations involving
point masses, point charges, and point vortices,
as well as in boundary element methods where the particles are quadrature points.
Evaluating~\eqref{N-body} by direct summation requires $O(N^2)$ operations
which is prohibitively expensive when $N$ is large,
and
several fast methods have been developed to reduce the cost.
One can distinguish between two types of methods,
{\it particle-mesh methods}
in which the particles are projected onto a uniform mesh where the FFT or multigrid can be used
(e.g.~P3M~\cite{HockneyEastwood},
particle-mesh Ewald~\cite{Essmann},
spectral Ewald~\cite{Klinteberg2017},
multilevel summation~\cite{brandt-lubrecht-90,hardy-skeel-16}),
and
{\it tree-based methods} 
in which the particles are partitioned into a hierarchy of clusters with a tree structure
and
the particle-particle interactions are replaced by
particle-cluster or cluster-cluster approximations 
(e.g.~treecode~\cite{BarnesHut}, 
fast multipole method (FMM)~\cite{GreengardRokhlin},
panel clustering~\cite{hackbusch-nowak-89}).

{\bf Tree-based methods.}
The present work is concerned with tree-based methods 
that rely on degenerate kernel approximations of the form,
\begin{equation}
\label{degenerate_kernel}
k({\bf x},{\bf y}) \approx \sum_{k=0}^n \phi_k({\bf x})\psi_k({\bf y}).
\end{equation}
Such approximations can be classified as 
{\it near-field/local} or {\it far-field/multipole} depending on their domain of validity
in the variables ${\bf x}, {\bf y}$.
The treecode originally used a far-field monopole approximation for the 
Newtonian potential~\cite{BarnesHut},
while the FMM improved on this by employing 
higher-order multipole and local approximations,
in particular using Laurent series for the 2D Laplace kernel
and
spherical harmonics for the 3D Laplace kernel~\cite{GreengardRokhlin,Greengard}.
Later versions of the FMM used 
plane wave expansions for the 3D Laplace kernel~\cite{cheng-greengard-rokhlin-99}
and
spherical Bessel function expansions for the Yukawa potential~\cite{GreengardHuang}.
Methods based on Cartesian Taylor expansions were also
developed for some common kernels~\cite{Draghicescu,duan-k-01,lindsay-krasny-01,tausch-03,LiJohnstonKrasny,KrasnyWang,WangTlupovaKrasny}.

{\bf Kernel-independent methods.}
The tree-based methods cited above rely on analytic series expansions specific to each kernel
and
alternative approximation methods have been investigated.
An early example in this direction was an FMM for Laplace kernels based on
discretizing the Poisson integral formula~\cite{anderson-92},
and
this was followed by a pseudoparticle method that reproduces the multipole moments 
for these kernels~\cite{makino-99}.
Later work developed approximations
suitable for a wide class of non-oscillatory kernels.
One approach based on polynomial interpolation~(\cite{kress-14}, section 11.4)
has been applied in the context of
multilevel approximation~\cite{giebermann-01},
hierarchical matrices~\cite{borm-grasedyck-hackbusch-03},
and
the black-box FMM (bbFMM)~\cite{FongDarve}.
An alternative method employed in the
kernel-independent FMM (KIFMM) uses equivalent densities~\cite{YingBirosZorin,Ying},
and
other kernel-independent FMMs use Legendre expansions~\cite{GimbutasRokhlin},
matrix compression based on skeletonization~\cite{martinsson-rokhlin-07},
and
truncated Fourier series~\cite{zhang-sun-11}.
Recently 
an FMM based on the Cauchy integral formula and Laplace transform was proposed 
for general analytic functions~\cite{letourneau-cecka-darve-14},
and
a kernel-independent treecode was developed using approximate skeletonization
for particle systems in high dimensions~\cite{march-xiao-biros-15}.

{\bf Present work.}
There is ongoing interest in exploring different strategies for fast summation of
particle interactions,
and
the present work contributes a kernel-independent treecode (KITC)
with operation count $O(N\log N)$
in which the far-field approximation uses 
barycentric Lagrange interpolation at Chebyshev points~\cite{BerrutTrefethen,trefethen-13}. 
The barycentric Lagrange interpolant can be efficiently implemented 
and
has good stability properties~\cite{rack-reimer-82,Higham,mascarenhas-14};
the 1D case is reviewed in~\cite{BerrutTrefethen,trefethen-13}
and
here we apply it in 3D using a tensor product to compute 
well-separated particle-cluster approximations.

It should be noted that the bbFMM~\cite{FongDarve} 
and 
KITC both use polynomial interpolation,
but they differ in two ways.
The first difference concerns the interpolating polynomial;
the bbFMM uses Chebyshev points of the 1st kind (roots of Chebyshev polynomials) 
and
expresses the Lagrange polynomials in terms of Chebyshev polynomials,
while the KITC uses Chebyshev points of the 2nd kind (extrema of Chebyshev polynomials)
and
expresses the Lagrange polynomials in barycentric form;
as explained in~\Cref{sec:barycentric},
this enables the KITC to take advantage of the scale-invariance property 
of barycentric Lagrange interpolation.
The second difference concerns the algorithm structure;
the KITC uses only far-field approximations
and
avoids the multipole-to-local translations
and
SVD compression steps in the bbFMM~\cite{FongDarve}.
With these choices the KITC
is a relatively simple algorithm with low memory consumption,
and
this enables a straightforward OpenMP parallelization.

We present numerical results motivated by the method of regularized Stokeslets (MRS)
for slow viscous flow~\cite{Cortez,CortezFauciMedovikov}. 
The MRS has been applied to simulate
cilia- and flagella-driven flow~\cite{Flores,Smith}, 
helical swimming~\cite{CortezFauciMedovikov}, 
slender body flow~\cite{BouzarthMinion},
coupled Stokes-Darcy flow~\cite{TlupovaCortez},
and
flow around elastic rods~\cite{OlsonLimCortez}.
Due to the complexity of the MRS kernels,
they are prime candidates for kernel-independent fast summation methods,
but as far as we know only recently has the KIFMM been applied to MRS simulations~\cite{RostamiOlson}.
Here we apply the KITC to systems of regularized Stokeslets and rotlets
from~\cite{RostamiOlson}.
The results demonstrate the method's good performance in terms of
accuracy, efficiency, and memory consumption in serial and parallel simulations.

The paper is organized as follows. 
\Cref{sec:polynomial_interpolation} discusses polynomial interpolation 
and 
its application to kernel approximation.
\Cref{sec:barycentric} reviews barycentric Lagrange interpolation following~\cite{BerrutTrefethen,trefethen-13}.
\Cref{sec:particle-cluster} explains how the interpolant is used to approximate 
particle-cluster interactions.
\Cref{sec:modified_weights} describes how some quantities called modified weights
are computed.
\Cref{sec:treecode} presents the KITC algorithm. 
\Cref{sec:MRS_kernels} reviews the MRS kernels (regularized Stokeslet and rotlet).
\Cref{sec:results} presents numerical results for two examples
motivated by recent MRS simulations~\cite{RostamiOlson}.
A summary is given in~\Cref{sec:summary}. 


\section{Polynomial interpolation and kernel approximation}
\label{sec:polynomial_interpolation}

We begin by recalling
some basic facts about polynomial interpolation in 1D~\cite{trefethen-13}.
Given a function $f(t)$
and
$n+1$ distinct points $s_k \in [-1,1]$ for $k = 0,\ldots,n$,
there is a unique polynomial $p_n(t)$ of degree at most $n$ that interpolates the function at these points,
$p_n(s_k) = f(s_k), k = 0,\ldots,n$.
The Lagrange form of the interpolating polynomial is
\begin{equation}
\label{Lagrange_form}
p_n(t) = \sum_{k=0}^n f_k L_k(t), \quad f_k = f(s_k), \quad k = 0,\ldots,n,
\end{equation}
where the Lagrange polynomials,
\begin{equation}
\label{Lagrange_polynomial}
L_k(t) = \frac{\Pi_{j=0,j \neq k}^n (t-s_j)}{\Pi_{j=0,j \neq k}^n (s_k-s_j)}, \quad k=0,\ldots,n,
\end{equation}
have degree $n$ 
and
satisfy $L_k(s_j) = \delta_{jk}$.
We view the interpolating polynomial $p_n(t)$ as an approximation to $f(t)$,
and
applying this idea to a kernel $k(x,y)$ in 1D,
we hold $x$ fixed and interpolate with respect to $y$ to obtain 
\begin{equation}
\label{separated_1D}
k(x,y) \approx \sum_{k=0}^n k(x,s_k) L_k(y).
\end{equation}
The approximation on the right is a polynomial of degree $n$ in the variable $y$
and
it interpolates the kernel at $y=s_k$;
this idea is well known (e.g.~\cite{kress-14}, section 11.4).

Now consider a kernel $k({\bf x},{\bf y})$ in 3D
and
a tensor product set of grid points ${\bf s}_{\bf k} = (s_{k_1}, s_{k_2}, s_{k_3}) \in [-1,1]^3$,
where ${\bf k} = (k_1,k_2,k_3)$ is a multi-index with $k_\ell = 0,\ldots,n$ for $\ell = 1,2,3$.
As above we hold ${\bf x}$ fixed and interpolate with respect to ${\bf y}=(y_1,y_2,y_3)$ to obtain
\begin{equation}
\label{separated_3D}
k({\bf x},{\bf y}) \approx 
\sum_{k_1=0}^n\sum_{k_2=0}^n\sum_{k_3=0}^n k({\bf x},{\bf s}_{\bf k}) 
L_{k_1}(y_1)L_{k_2}(y_2)L_{k_3}(y_3).
\end{equation}
In this case the approximation on the right 
is a polynomial of degree $n$ in each variable $(y_1,y_2,y_3)$
and
it interpolates the kernel at the grid points ${\bf y} = {\bf s}_{\bf k}$.
Moreover, the Lagrange polynomial expressions~\eqref{separated_1D} and \eqref{separated_3D}
are degenerate kernel approximations of the form~\eqref{degenerate_kernel}.

\section{Barycentric Lagrange interpolation}
\label{sec:barycentric}

The expression for the interpolating polynomial described
in~\Cref{sec:polynomial_interpolation}
is not well-suited for practical computing due to cost and stability issues;
the problem though is not with the Lagrange form~\eqref{Lagrange_form} for $p_n(t)$,
but rather with the expression~\eqref{Lagrange_polynomial} for the Lagrange polynomials~$L_k(t)$.
Berrut and Trefethen~\cite{BerrutTrefethen} advocated using instead the 
2nd barycentric form of the Lagrange polynomials,
\begin{equation}
\label{Barycentric-1D}
L_k(t) = \frac{\displaystyle \frac{w_k}{t - s_k}}
{\sum\limits_{k = 0}^n \displaystyle \frac{w_k}{t - s_k}}, \quad
w_k = \frac{1}{\Pi_{j=0,j\neq k}^n (s_k - s_j)}, \quad k = 0,\ldots,n,
\end{equation}
where $w_k$ are the barycentric weights.
This form is mathematically equivalent to~\eqref{Lagrange_polynomial},
with the understanding that the removable singularity at $t=s_j$ 
is resolved by setting $L_k(s_j) = \delta_{jk}$;
to enforce this condition,
following~\cite{BerrutTrefethen} the code is written so that
if the argument $t$ is closer to an interpolation point $s_j$ than some tolerance,
this is flagged 
and 
the correct value $L_k(s_j) = \delta_{jk}$ is provided.
In our implementation the tolerance is the  
minimum positive IEEE double precision floating point number,
DBL\_MIN = 2.22507e-308;
further details will be given in~\Cref{sec:modified_weights}.

We work with Chebyshev points of the 2nd kind,
\begin{equation}
\label{Chebyshev_points}
s_k = \cos\theta_k, \quad \theta_k = \frac{k\pi}{n},\quad  k = 0,\ldots,n.
\end{equation}
In this case the interpolating polynomial $p_n(t)$ converges rapidly 
and uniformly on $[-1,1]$ to $f(t)$ as $n$ increases,
under mild smoothness assumptions on the given function~\cite{trefethen-13}.
Note that computing the barycentric weights~$w_k$ 
by the definition~\eqref{Barycentric-1D} requires $O(n^2)$ operations,
but this expense disappears for the $s_k$ chosen in~\eqref{Chebyshev_points}
because in that case the following simple weights can be 
used instead~\cite{Salzer,trefethen-13},
\begin{equation}
\label{barycentric_weights}
w_k = (-1)^k\delta_k, \quad \delta_k = 
\begin{cases} 
1/2, & k = 0 ~\text{or}~ k=n, \cr
1, & k = 1,\ldots,n-1.
\end{cases}
\end{equation}
This relies on a scale-invariance property of the barycentric form of $L_k(t)$ 
in~\eqref{Barycentric-1D};
namely, if the weights $w_k$ have a common constant factor $\alpha \ne 0$,
then $\alpha$ can be cancelled from the numerator and denominator,
and~\eqref{Barycentric-1D} stays the same.

Hence the barycentric Lagrange form of the interpolating polynomial is
\begin{equation}
p_n(t) = \sum_{k=0}^n \frac
{\displaystyle \frac{w_k}{t - s_k}}
{\displaystyle \sum_{k=0}^n\frac{w_k}{t - s_k}} f_k,
\label{barycentric_Lagrange_form}
\end{equation}
and
with the simple weights~\eqref{barycentric_weights}, 
evaluating $p_n(t)$ requires $O(n)$ operations~\cite{BerrutTrefethen,trefethen-13}.

The scale-invariance property is also important when working on different intervals;
if $[-1,1]$ is linearly mapped to $[a,b]$ by $t \to \frac{1}{2}(a+b + t(b-a))$,
then the weights $w_k$ defined in~\eqref{Barycentric-1D} 
gain a factor of $2^n/(b-a)^n$ which could lead to overflow or underflow,
but this factor can be safely omitted due to scale-invariance.
This means that the simple weights~\eqref{barycentric_weights}
can be used for any interval $[a,b]$,
along with the linearly mapped Chebyshev points;
this is important in the present work because the treecode 
uses intervals of different sizes.
Note also that barycentric Lagrange interpolation
is stable in finite precision 
arithmetic~\cite{rack-reimer-82,Higham,mascarenhas-14},
and
the Chebfun software package uses this form of
polynomial interpolation~\cite{driscoll-hale-trefethen-14}.

For comparison,
the bbFMM~\cite{FongDarve} uses Chebyshev points of the 1st kind,
\begin{equation}
\bar{s}_k = \cos\bar{\theta}_k, \quad \bar{\theta}_k = \frac{(2k-1)\pi}{2n}, \quad k = 1,\ldots,n,
\end{equation}
and 
a Chebyshev Lagrange form of the interpolating polynomial,
\begin{equation}
\label{Chebyshev}
p_{n-1}(t) = \sum_{k=1}^n f_k \bar{L}_k(t), \quad
\bar{L}_k(t) = \frac{1}{n} + \frac{2}{n}\sum_{j=1}^{n-1}T_j(\bar{s}_k)T_j(t),
\quad k = 1,\ldots,n,
\end{equation}
where $T_j(t)$ is the $j$th degree Chebyshev polynomial
and
$\bar{L}_k(\bar{s}_j) = \delta_{jk}$.
The cost of evaluating $p_{n-1}(t)$ by 
directly summing~\eqref{Chebyshev} is $O(n^2)$,
or $O(n\log n)$ if a fast transform is used.
We are not aware of an analog of the barycentric scale-invariance property.



\section{Particle-cluster interactions}
\label{sec:particle-cluster}

\begin{figure}[htb]
\centering
\includegraphics[width=0.95\textwidth]{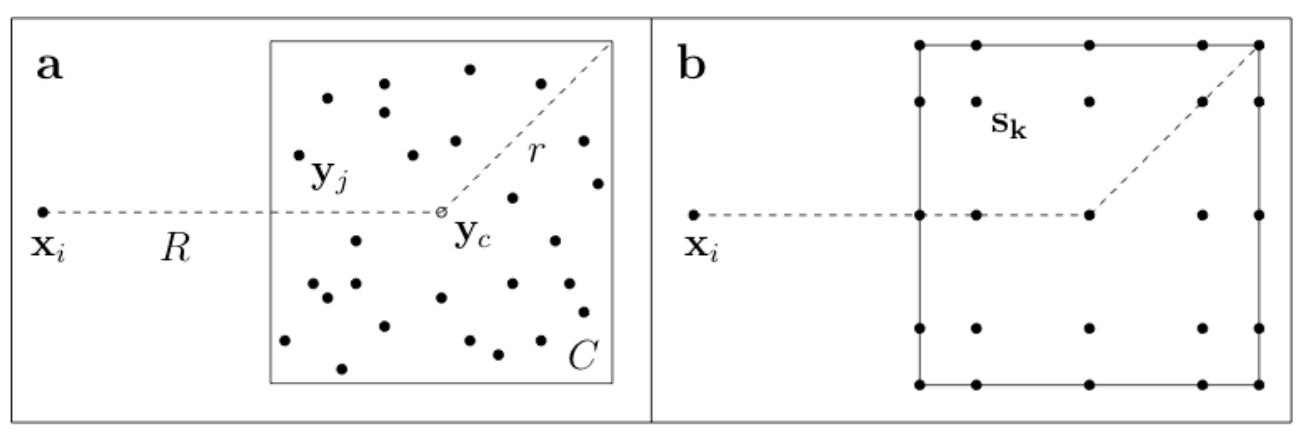}
\caption{
Particle-cluster interaction,
(a) target particle ${\bf x}_i$,
source cluster $C = \{{\bf y}_j\}$,
center ${\bf y}_c$, 
radius $r$, 
particle-cluster distance $R = |{\bf x}_i - {\bf y}_c|$,
(b) Chebyshev grid points ${\bf s}_{\bf k}$ in cluster $C$.}
\label{fig:far_field_schematic}
\end{figure}

In a treecode the particles ${\{\bf x}_i\}$ are partitioned into a hierarchy of clusters $\{C\}$,
and
the sum~\eqref{N-body} is written as
\begin{equation}
\label{N-body_2}
u({\bf x}_i) = 
\sum_{j=1}^N k({\bf x}_i,{\bf x}_j) f_j =
\sum_C u({\bf x}_i,C), \quad i = 1,\ldots,N,
\end{equation}
where
\begin{equation}
u({\bf x}_i,C) = \sum_{{\bf y}_j \in C} k({\bf x}_i,{\bf y}_j) f_j,
\label{pc_interaction_3D}
\end{equation}
is the interaction between a target particle ${\bf x}_i$
and
a source cluster $C = \{{\bf y}_j\}$.
The sum over $C$ in~\eqref{N-body_2} denotes a suitable subset of clusters
depending on the target particle.
Figure~\ref{fig:far_field_schematic}a depicts a particle-cluster interaction with
target particle~${\bf x}_i$,
cluster center~${\bf y}_c$,
cluster radius~$r$,
and
particle-cluster distance $R = |{\bf x}_i - {\bf y}_c|$.
In this work the clusters are rectangular boxes whose sides are aligned with the 
coordinate axes.
When the particle and cluster are well-separated (the criterion is given in \Cref{sec:treecode}), 
the interaction~\eqref{pc_interaction_3D} is computed using the 
kernel approximation~\eqref{separated_3D};
this is depicted in Fig.~\ref{fig:far_field_schematic}b where the Chebyshev grid points 
${\bf s}_{\bf k}$ are mapped to the cluster.
Hence the target particle~${\bf x}_i$ interacts with the interpolation points ${\bf s}_{\bf k}$
rather than the source particless~${\bf y}_j$.


In more detail,
we substitute the kernel approximation~\eqref{separated_3D}
into the particle-cluster interaction~\eqref{pc_interaction_3D}
and
switch the order of summation
to obtain the far-field particle-cluster approximation,
\begin{equation}
\label{pc_approximation_3D}
u({\bf x}_i,C) \approx
\sum_{k_1=0}^n \sum_{k_2=0}^n \sum_{k_3=0}^n k({\bf x}_i, {\bf s}_{\bf k}) \widehat{f}_{\bf k},
\end{equation}
where the modified weights are
\begin{equation}
\label{modified_weights}
\widehat{f}_{\bf k} =  
\sum_{{\bf y}_j \in C}
L_{k_1}(y_{j1}) L_{k_2}(y_{j2}) L_{k_3}(y_{j3})f_j, 
\quad k_1,k_2,k_3 = 0,\ldots,n, 
\end{equation}
and ${\bf y}_j = (y_{j1},y_{j2},y_{j3})$.
Note that when the target ${\bf x}_i$ is well-separated from the 
sources ${\bf y}_j$,
the kernel $k({\bf x},{\bf y})$ is interpolated on a subdomain where it is smooth 
and
this ensures the accuracy of the approximation.

We see that~\eqref{pc_approximation_3D} has the form of a 
particle-cluster interaction~\eqref{pc_interaction_3D},
where the source particles and weights $\{{\bf y}_j,f_j\}$ 
are replaced by the Chebyshev grid points and modified weights $\{{\bf s}_k,\widehat{f}_{\bf k}\}$. 
This is advantageous for two reasons.
First, the modified weights $\widehat{f}_{\bf k}$ for a given cluster 
are independent of the target particle ${\bf x}_i$, 
so they can be precomputed and re-used for different targets;
details are given below.
Second, the cost of the direct sum~\eqref{pc_interaction_3D} is $O(N_c)$, 
where $N_c$ is the number of particles in cluster $C$, 
while the cost of the approximation~\eqref{pc_approximation_3D} is $O(n^3)$ 
assuming the modified weights are known,
so there is a cost reduction when $N_c\,{>>}\,n^3$.
Finally note that the approximation~\eqref{pc_approximation_3D} 
depends only on kernel evaluations
and
hence qualifies as kernel-independent.

\section{Computation of modified weights}
\label{sec:modified_weights}
Algorithm~\ref{modified_weights_algorithm} describes how the
modified weights are computed.
Before proceeding,
note that the Lagrange polynomial terms in~\eqref{modified_weights} 
in barycentric form are
\begin{equation}
L_{k_\ell}(y_{j\ell}) = 
\frac{\displaystyle \frac{w_{k_\ell}}{y_{j\ell} - s_{k_\ell}}}
{\sum\limits_{k = 0}^n \displaystyle \frac{w_k}{y_{j\ell} - s_k}} =
\frac{a_{j,k_\ell,\ell}}{\displaystyle \sum_{k=0}^n a_{j,k,\ell}}, \quad
a_{j,k,l} = \frac{w_k}{y_{jl} - s_k},
\end{equation}
where the variables $a_{j,k,l}$ are introduced for clarity.
There are three loops,
the outer loop over source particles $j = 1,\ldots,N_c$,
the 1st inner loop over Chebyshev points $k = 0,\ldots,n$,
and
the 2nd inner loop over coordinate indices $\ell = 1,2,3$.
To handle the removable singularity,
we adapt the procedure in~\cite{BerrutTrefethen} (page 510).

In Algorithm~\ref{modified_weights_algorithm},
the source particles and weights ${\bf y}_j, f_j$ for a given cluster are input,
and
the modified weights $\widehat{f}_{\bf k}$ are output.
The modified weights are initialized to zero on line 4,
the loop over source particles starts on line 6,
and
the flags and sums are initialized on line 7.
Lines 9-13 loop over the Chebyshev points and coordinate indices to compute the terms 
$a_{j,k,\ell}$;
this requires $O(n N_c)$ operations
and
a temporary storage array of size $O(nN_c)$ which is reused for different clusters.
Line 10 checks whether a source particle coordinate is close to a Chebyshev point;
if so, then a flag records the index.
Then on lines 15-17 the code checks to see whether a flag was set;
if so, then the temporary variables are adjusted to enforce the 
condition $L_k(s_j) = \delta_{jk}$.
Then on lines 20-22 the code loops over the tensor product of Chebyshev point indices
and
increments the modified weights~$\widehat{f}_{\bf k}$ as indicated in~\eqref{modified_weights};
this requires $O(n^3 N_c)$ operations
and $O(n^3)$ storage for each cluster.
The modified weights are computed this way for each cluster when the tree is constructed
and
in practice this requires only a small fraction of the total KITC CPU time;
for example in one case with $N$ = 640\,K particles, 
the total CPU time was 149\,s, 
while computing the modified weights took less than 2\,s.

\begin{algorithm}[htb]
\caption{computation of modified weights in~\eqref{modified_weights}}
\label{modified_weights_algorithm}
\begin{algorithmic}[1]
\State input: source particles and weights for a given cluster, ${\bf y}_j, f_j$ 
\State input: Chebyshev points mapped to the cluster, $s_k$ 
\State output: modified weights, $\widehat{f}_{\bf k}$
\State initialize all $\widehat{f}_{\bf k} = \widehat{f}(k_1,k_2,k_3) = 0$
\State \% loop over source particles
\State for $j = 1:N_c$
\State \quad initialize flag(1:3) = -1, sum(1:3) = 0
\State \% loop over Chebyshev points and coordinate indices
\State \quad for $k = 0:n$ and $\ell = 1:3$
\State \qquad if $|y_{j\ell} - s_k| \le$ DBL\_MIN, flag$(\ell)$= $k$
\State \qquad\quad else $a(j,k,\ell) = w_k /(y_{j\ell} - s_k)$, sum$(\ell)$ += $a(j,k,\ell)$
\State \qquad end if
\State \quad end for
\State \% if a flag was set, adjust sum$(\ell)$ and $a(j,k,\ell)$ to handle removable singularity
\State \quad for $\ell = 1:3$
\State \qquad if flag$(\ell) > -1$, sum$(\ell) = 1, a(j,0:n,\ell) = 0, a(j,\,$flag$(\ell), \ell) = 1$, end if
\State \quad end for
\State \quad denom = sum(1) $\cdot$ sum(2) $\cdot$ sum(3)
\State \% loop over tensor product of Chebyshev point indices as in~\eqref{modified_weights}
\State \quad for $(k_1,k_2,k_3) = (0:n,0:n,0:n)$
\State \qquad $\widehat{f}(k_1,k_2,k_3)$ += 
$(a(j,k_1,1) \cdot a(j,k_2,2) \cdot a(j,k_3,3) / {\rm denom}) \cdot f_j$
\State \quad end for
\State end for
\end{algorithmic}
\end{algorithm}

The quantities $L_{k_\ell}(y_{j\ell})$ could be computed instead using the 
Chebyshev form~\eqref{Chebyshev},
but the argument of the Chebyshev polynomials appearing there 
lies in the unit interval $[-1,1]$,
so the particles ${\bf y}_j$ must be mapped to the unit cube $[-1,1]^3$,
with cost $O(N_c)$.
The barycentric form~\eqref{Barycentric-1D} does not require this mapping;
the particles stay in place in the cluster;
hence the barycentric form has an advantage over the Chebyshev form
in that it avoids the cost of mapping the particles ${\bf y}_j$.

%


\section{Kernel-independent treecode algorithm}
\label{sec:treecode}

Aside from using barycentric Lagrange interpolation for the far-field approximation,
the present algorithm is similar to previous treecodes based on
analytic series expansions~\cite{BarnesHut,Draghicescu,LiJohnstonKrasny}.
The procedure is outlined in Algorithm~\ref{treecode_algorithm}. 
After inputting the particle data and treecode parameters,
a hierarchical tree of particle clusters is built
and
the modified weights for each cluster are computed.
In our implementation
the root cluster of the tree is the smallest rectangular box enclosing the particles.
The root is bisected in any coordinate direction for which the side length
is greater than $\ell_{max}/\sqrt{2}$,
where $\ell_{max}$ is the maximum side length of the root.
Hence the root is bisected in its long directions, 
but not its short directions, and the child clusters are bisected in the same way.
The process continues until a cluster has fewer than $N_0$ particles, 
a user-specified parameter.
The clusters obtained this way are rectangular boxes
and
a cluster may have 8, 4, or 2 children depending on which sides were bisected.
The flexibility to use rectangular boxes instead of cubes is important in
cases like Example 2 below where the particles lie in a rectangular slab.

\begin{algorithm}
\caption{kernel-independent treecode}
\label{treecode_algorithm}
\begin{algorithmic}[1]
\State input: particle coordinates and weights ${\bf x}_i, f_i, i=1,\ldots,N$ 
\State input: treecode MAC parameter $\theta$, polynomial degree $n$, maximum leaf size $N_0$
\State output: particle velocities $u_i, i=1,\ldots,N$
\State program {\bf main}
\State \quad build tree of particle clusters
\State \quad compute modified weights $\widehat{f}_{\bf k}$ in \eqref{modified_weights} for each cluster
\State \quad for $i = 1,\ldots,N$, {\bf compute\_velocity}(${\bf  x}_i$, root), end for
\State end program
\State subroutine {\bf compute\_velocity}(${\bf x}$, $C$)
\State \quad if MAC is satisfied
\State \qquad compute particle-cluster interaction by approximation~\eqref{pc_approximation_3D}
\State \quad else
\State \qquad if $C$ is a leaf, compute particle-cluster interaction by direct sum~\eqref{pc_interaction_3D}
\State \quad else
\State \qquad for each child $C^\prime$ of $C$, {\bf compute\_velocity}($\bf x$, $C^\prime$), end for
\State end subroutine
\end{algorithmic}
\end{algorithm}


Turning to Algorithm~\ref{treecode_algorithm},
after building the tree of clusters (line 5),
the code cycles through the particles (line 7),
and each particle interacts with clusters, 
starting at the root and proceeding to the child clusters.
The multipole acceptance criterion (MAC),
\begin{equation}
\frac{r}{R} \leq \theta,
\label{MAC}
\end{equation}
determines whether a particle ${\bf x}_i$ and cluster $C$ are well-separated,
where (recall Fig.~\ref{fig:far_field_schematic})
$r$ is the cluster radius,
$R$ is the particle-cluster distance,
and
$\theta$ is a user-specified parameter.
If the MAC (\ref{MAC}) is satisfied, the particle-cluster interaction is computed 
by barycentric Lagrange interpolation~\eqref{pc_approximation_3D} 
with user-specified degree $n$;
otherwise the code checks the child clusters, 
or if the cluster is a leaf (no children), 
then the interaction is computed directly by~\eqref{pc_interaction_3D}.

This is essentially the Barnes-Hut algorithm~\cite{BarnesHut}, 
extended to higher-order particle-cluster approximations computed by 
barycentric Lagrange interpolation.
There are two stages;
stage~1 is building the tree
and computing the modified weights for each cluster,
and
stage~2 is computing the particle-cluster interactions.
In principle both stages scale like $O(N\log N)$,
although in practice stage~1 is much faster than stage~2.
There are three user-specified parameters ($\theta, n, N_0$) that control
the accuracy and efficiency of the treecode;
the present work uses representative values with no claim that they are optimal.


\section{Regularized Stokes kernels}
\label{sec:MRS_kernels}

The method of regularized Stokeslets (MRS) uses regularized kernels to represent 
point forces and torques in slow viscous flow~\cite{Cortez,CortezFauciMedovikov,OlsonLimCortez}.
In particular, consider the following functions defined on~${\mathbb R}^3$,
\begin{subequations}
\label{MRS_kernels}
\begin{align}
H_1(r) &= \frac{2\epsilon^2 + r^2}{8\pi(r^2 + \epsilon^2)^{3/2}}, \quad 
H_2(r) = \frac{1}{8\pi(r^2 + \epsilon^2)^{3/2}}, \quad
Q(r) = \frac{5\epsilon^2 + 2r^2}{8\pi(r^2 + \epsilon^2)^{5/2}}, \\
D_1(r) &= \frac{10\epsilon^4 - 7\epsilon^2r^2 - 2r^4}{8\pi(r^2 + \epsilon^2)^{7/2}}, \quad 
D_2(r) = \frac{21\epsilon^2 + 6r^2}{8\pi(r^2 + \epsilon^2)^{7/2}},
\end{align}
\end{subequations}
where 
$r = |{\bf x}-{\bf y}|$,
${\bf x}$ is a target, 
${\bf y}$ is a source,
and
$\epsilon$ is the MRS regularization parameter.
Given a set of source particles $\{{\bf y}_j\}$
with forces $\{{\bf f}_j\}$ and torques $\{{\bf n}_j\}$, for $j = 1,\ldots,N$,
the velocity induced by regularized Stokeslets is
\begin{equation}
{\bf u}({\bf x}) = \sum_{j = 1}^N
\left({\bf f}_jH_1(r_j) + \left[{\bf f}_j\cdot ({\bf x} - {\bf y}_j)\right] ({\bf x} - {\bf y}_j)H_2(r_j)\right),
\label{Velocity_1}
\end{equation}
where $r_j = |{\bf x} - {\bf y}_j|$,
while the linear velocity and angular velocity induced by
regularized Stokeslets and rotlets are
\begin{subequations}
\label{MRS_2}
\begin{align}
{\bf u}({\bf x}) &= \sum_{j = 1}^N
\bigg({\bf f}_jH_1(r_j) + \left[{\bf f}_j\cdot ({\bf x} - {\bf y}_j)\right] ({\bf x} - {\bf y}_j)H_2(r_j) + 
\frac{1}{2}\left[{\bf n}_j\times({\bf x} - {\bf y}_j)\right]Q(r_j)\bigg),
\label{Velocity_2} \\
{\bf w}({\bf x}) &= \sum_{j = 1}^N
\bigg( \frac{1}{2}\left[{\bf f}_j\times({\bf x} - {\bf y}_j)\right]Q(r_j) + 
\frac{1}{4}{\bf n}_jD_1(r_j) + \frac{1}{4}\left[{\bf n}_j\cdot({\bf x} - {\bf y}_j) \right] D_2(r_j)\bigg).\label{AngularVelocity}
\end{align}
\end{subequations}

These sums have the pairwise interaction form~\eqref{N-body}
suitable for a fast summation method;
in this case the kernels are tensors,
and
the particle weights and output velocities are vectors.
However the MRS functions~\eqref{MRS_kernels} are somewhat complicated
and
we know of only one application~\cite{RostamiOlson}
using a version of the KIFMM in which the equivalent densities
are defined on coronas (or shells) around each cluster~\cite{Ying}.
Note also that when $\epsilon\,{=}\,0$,
the MRS kernels reduce to the usual Stokes kernels which are homogeneous
(i.e.~$k(\alpha{\bf x},\alpha{\bf y}) = \alpha^\lambda k({\bf x},{\bf y})$ for some constant $\lambda$
and
all $\alpha > 0$); 
some versions of the FMM use this property to improve performance 
(e.g.~\cite{malhotra-biros-16,takahashi-cecka-darve-12}),
but this optimization is not available for the MRS kernels 
because they are non-homogeneous when $\epsilon \ne 0$.
In the next section we demonstrate the capability of the KITC in evaluating the 
MRS sums~\eqref{Velocity_1}-\eqref{MRS_2} with $\epsilon \ne 0$.


\section{Numerical results}
\label{sec:results}

We present results for two examples motivated by recent MRS simulations~\cite{RostamiOlson}.
In these examples the targets and sources coincide,
but this is not an essential restriction.
To quantify the accuracy of the KITC we define the relative error,
\begin{equation}
\label{error}
E = 
\left(\sum\limits_{i=1}^N |{\bf u}^d({\bf x}_i) - {\bf u}^t({\bf x}_i)|^2\bigg/
\sum\limits_{i=1}^N |{\bf u}^d({\bf x}_i)|^2\right)^{\!\!1/2},
\end{equation}
where ${\bf u}^d({\bf x}_i)$ is the exact velocity computed by direct summation,
${\bf u}^t({\bf x}_i)$ is the treecode approximation,
and $|{\bf u}|$ is the Euclidean norm.
In Example 2 where the linear velocity~\eqref{Velocity_2}
and 
angular velocity~\eqref{AngularVelocity} are computed,
they are combined into a single vector
and
the error is computed as in~\eqref{error}.
All lengths are nondimensional.
References to \lq\lq CPU time" are the total wall-clock run time in seconds.

The algorithm was coded in C++ in double precision
and
compiled using the Intel compiler icpc with $-$O2 optimization.
The computations were performed on the University of Wisconsin-Milwaukee 
Mortimer Faculty Research Cluster;
each node is a Dell PowerEdge R430 server with 
two 12-core Intel Xeon E5-2680 v3 processors at 2.50\,GHz and 64\,GB RAM;
all compute and I/O nodes are linked by Mellanox FDR Infiniband (56Gb/s) 
and gigabit Ethernet networks.
Parallel runs used OpenMP with up to 24 cores on a single node.
The source code is available for download (github.com/Treecodes/stokes-treecode).


\subsection{Example 1}

The first example simulates microorganisms randomly located
in a cube of side length $L=10$,
where each microorganism is a pair of particles representing its body and flagella,
and
the particles exert unit forces in opposite directions along the organism length~\cite{RostamiOlson}. 
The total number of particles is $N$
and
we consider five systems with $N\,{=}\,$10K, 80K, 640K, 5.12M, 40.96M.
The maximum number of particles in a leaf was set to $N_0 = 2000$;
note however that in this case each cluster has 8 children
and
the actual leaf size is 
approximately $N/8^d = 1250$, where $d$ is the depth of the tree
(some leaves may have a few more or less since the particle locations are random).
The microorganism length is $\ell\,{=}\,0.02$
and
the MRS parameter is $\epsilon = 0.02$.
In this example we compute the Stokeslet velocity~\eqref{Velocity_1}.


\subsubsection{CPU time versus error}

\begin{figure}[htb]
\centering
\includegraphics[width=0.7\textwidth]{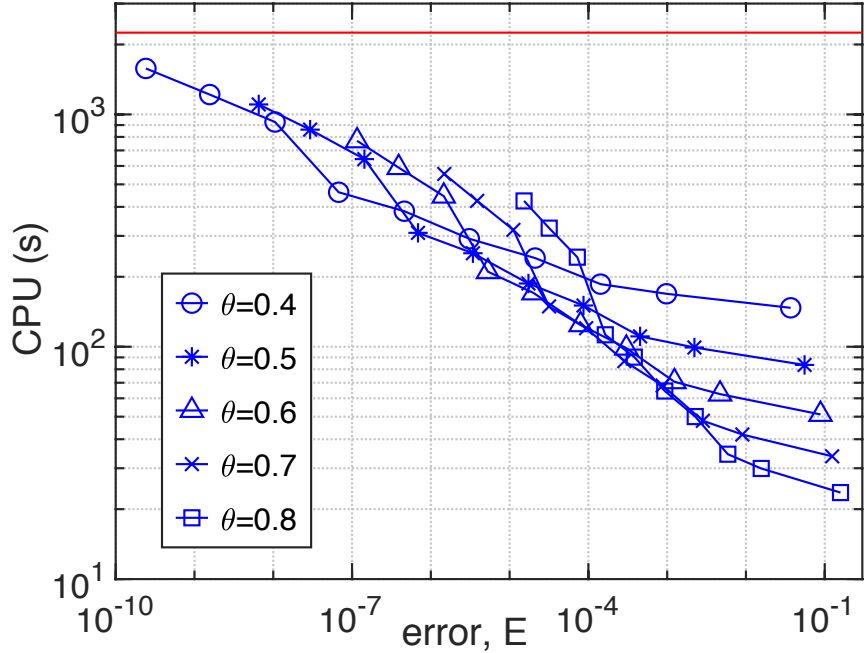}
\caption{Example 1,
regularized Stokeslets  in a cube~\eqref{Velocity_1},
MRS parameter $\epsilon = 0.02$,
system size $N\,{=}\,$640K,
direct sum CPU time (2249\,s, red horizontal line), 
treecode CPU time (s) is plotted versus error $E$ (blue symbols), 
MAC parameter $0.4 \le \theta \le 0.8$,
degree $n = 1,\ldots,10$ (increasing from right to left).}
\label{fig:CaseI_CPU_ERR_640K_DiffTheta}
\end{figure}

Figure~\ref{fig:CaseI_CPU_ERR_640K_DiffTheta} 
focuses on intermediate system size $N\,{=}\,$640K,
showing the treecode CPU time versus the error $E$,
for MAC parameters in the range $0.4 \le \theta \le 0.8$
and
interpolation degrees $n = 1,\ldots,10$ (increasing from right to left). 
For a given MAC parameter $\theta$,
the error $E$ decreases as the degree $n$ increases, 
but the CPU time increases. 
The lower envelope of the data gives the most efficient treecode parameters.
For example
to achieve error $E \le $\,1e-4, we can choose MAC $\theta = 0.7$ and degree $n=7$.
Figure~\ref{fig:CaseI_CPU_ERR_640K_DiffTheta} also shows the 
direct sum CPU time (2249\,s) as a horizontal line
(the value is independent of the error); 
the treecode is faster than direct summation for this range of parameters.

\subsubsection{Error and CPU time versus system size}

\begin{figure}[htb]
\centering
\includegraphics[width=0.495\textwidth]{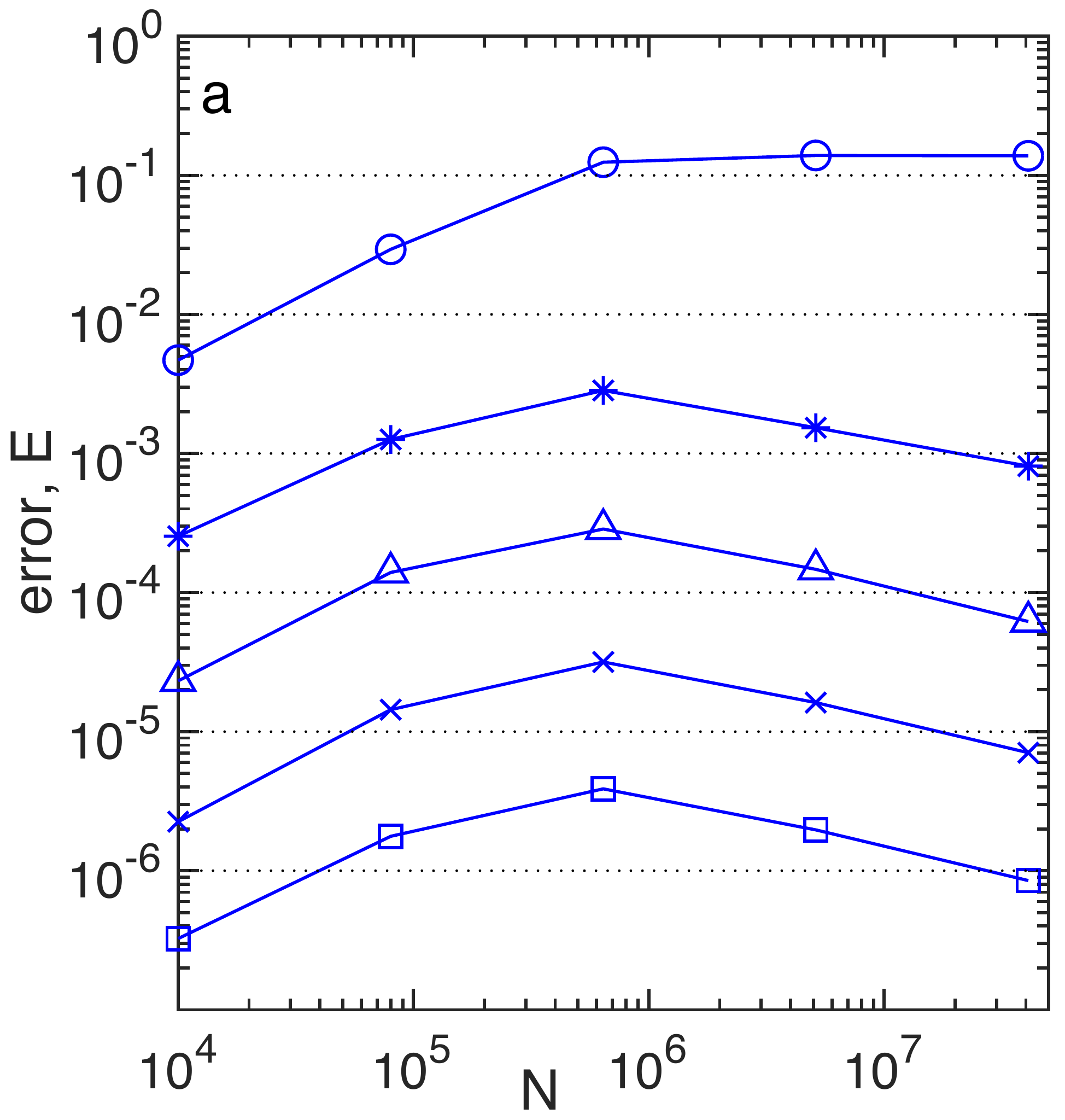}
\includegraphics[width=0.495\textwidth]{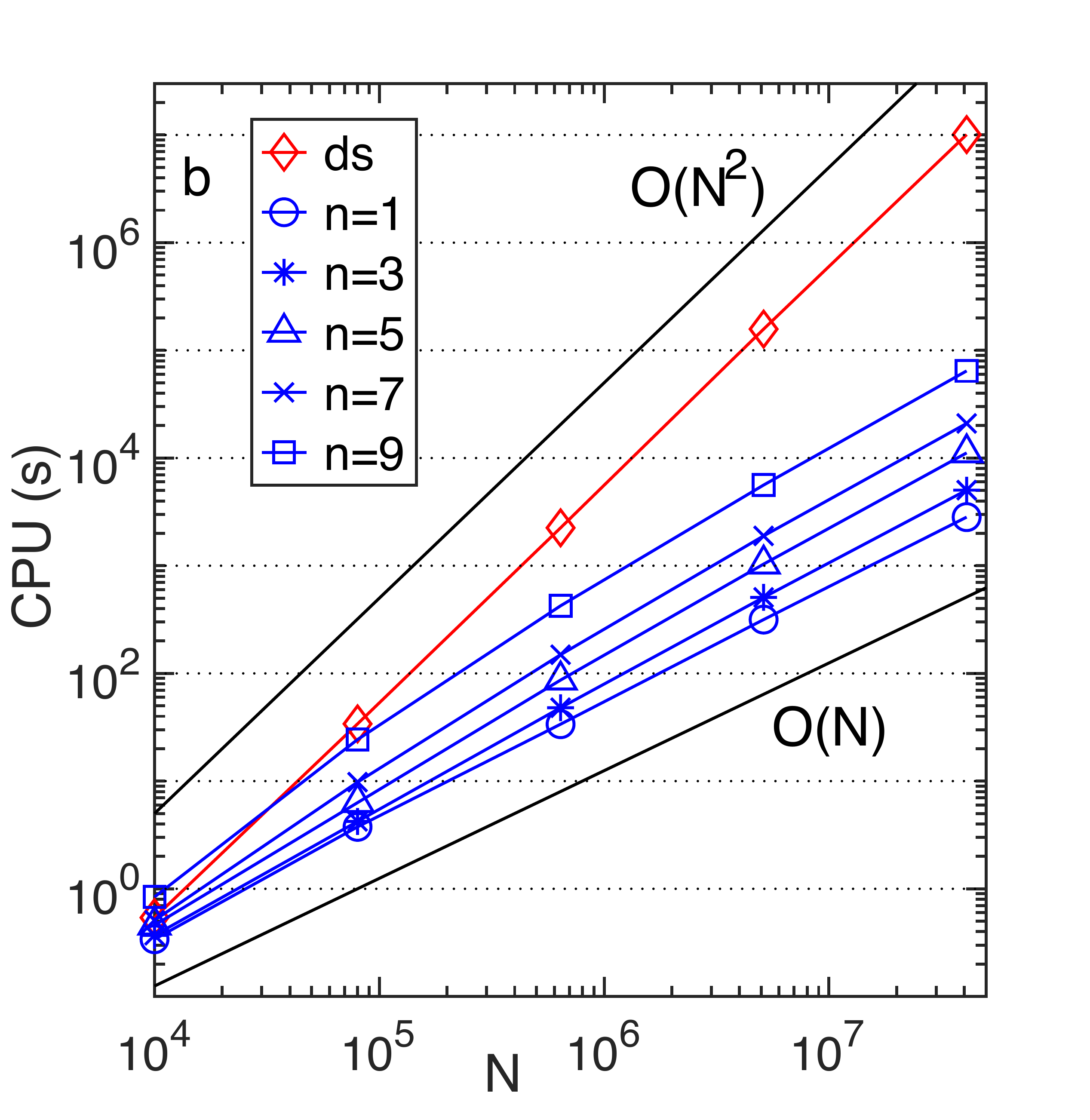}
\caption{Example 1,
regularized Stokeslets  in a cube~\eqref{Velocity_1},
MRS parameter $\epsilon = 0.02$,
(a) treecode error $E$,
(b) CPU time (s)
versus system size $N$ (10K, 80K, 640K, 5.12M, 40.96M),
treecode MAC parameter $\theta = 0.7$, degree $n = 1,3,5,7,9$, maximum leaf size $N_0=2000$,
direct sum CPU time (red,$\Diamond$).
}
\label{fig:CaseI_CPUErr_N}
\end{figure}

Figure~\ref{fig:CaseI_CPUErr_N} plots the treecode error $E$ (a),
and
CPU time (b) versus system size $N$,
for MAC parameter $\theta = 0.7$ and interpolation degree $n = 1,3,5,7,9$. 
Since it is impractical to carry out the direct sum for the largest system with $N\,{=}\,$40.96M,
in that case the error was computed using a sample of 2000 particles
and 
the CPU time was extrapolated from smaller systems.
In Fig.~\ref{fig:CaseI_CPUErr_N}a for a given degree $n$,
the error first increases with increasing system size~$N$
and
then it decreases slightly for the two largest systems;
nonetheless for a given system size $N$,
the error decreases as the interpolation degree $n$ increases.

Figure~\ref{fig:CaseI_CPUErr_N}b shows that the treecode is faster than direct summation 
except for the smallest system with $N\,{=}\,$10K
(those runs all take less than 1\,s).
Two reference lines are shown;
the direct sum CPU time is parallel to the slope 2 line consistent with $O(N^2)$ scaling,
while the KITC time is slightly steeper than the slope 1 line consistent with $O(N\log N)$ scaling.
As expected, for a given degree $n$
the speedup provided by the treecode increases with the system size~$N$.
To quantify this observation,
Table~\ref{table:CaseI_new} records the direct sum and treecode CPU time
and
the treecode error for the five systems using MAC parameter $\theta = 0.7$ and degree $n\,{=}\,7$.
For the largest system with $N\,{=}\,$40.96M particles,
the treecode is more than 480 times faster than direct summation while 
achieving error $E\,{=}\,$7.06e-6.

\begin{table}[htb]
\centering 
\begin{tabular}{|r|r|r|r|r|}
\hline
$N~~~$ & direct sum (s), $d$ & treecode (s), $t$ & speedup, $d/t$ & error, $E$ \\
\hline
10K        & 0.54~~                    & 0.51~~~              & 1.06~~~~          & 2.25e-6 \\
\hline
80K        & 34.01~~                  & 9.78~~~              & 3.48~~~~          & 1.44e-5 \\
\hline
640K      & 2248.93~~            & 149.14~~~          & 15.08~~~~        & 3.17e-5 \\
\hline
5.12M    & 157,542.13~~         & 1892.57~~~         & 83.24~~~~        & 1.62e-5 \\
\hline
40.96M  & 10,082,696.32~~    & 20,966.04~~~      & 480.91~~~~      & 7.06e-6 \\
\hline    
\end{tabular}
\caption{
Example 1,
regularized Stokeslets  in a cube~\eqref{Velocity_1},
MRS parameter $\epsilon = 0.02$,
system size $N$,
CPU time (s) for direct sum and treecode ($d,t$),
speedup ($d/t$),
treecode error ($E$),
results are from Fig.~\ref{fig:CaseI_CPUErr_N} with
MAC parameter $\theta = 0.7$, degree $n = 7$, maximum leaf size $N_0=2000$.
}
\label{table:CaseI_new}
\end{table}

 
\subsubsection{Memory consumption}  

The treecode memory consumption has two main components which are estimated as follows.
First, the memory due to the particle coordinates and weights is $O(N)$,
where $N$ is the number of particles in the system;
this is the same as in direct summation.
Second, the memory due to the modified weights is $O(n^3N)$,
where $n$ is the degree of polynomial interpolation
and
in this case $N$ represents the number of clusters in the tree
(this assumes the number of levels in the tree is $O(\log N)$).
Hence the theoretical memory consumption of the treecode is $O(n^3N)$,
although in practice the $O(N)$ particle data dominates when $n$ is small.

To see the actual values,
Table~\ref{tab:TreecodePerformance_Memory} 
presents the direct sum and treecode memory consumption (MB) for Example 1;
the results were obtained using the Valgrind Massif tool (www.valgrind.org).
For the largest system with $N\,{=}\,$40.96M particles,
the direct sum memory consumption (3604.80MB) was obtained by extrapolating
from the smaller systems.
It can be verified that the memory consumption in Table~\ref{tab:TreecodePerformance_Memory}
is consistent with the estimates in the previous paragraph;
for a given degree $n$, the memory consumption scales like $O(N)$,
while for a given system size $N$,
beyond the baseline particle data,
the additional memory consumption scales like $O(n^3)$.
Over this parameter range, 
the treecode uses less than 1.65 times as much memory as direct summation,
and
a KITC computation with $N\,{=}\,$40.96M particles and degree $n=9$
uses about 5.2GB of memory. 
For comparison,
a bbFMM computation with 1M~Stokeslets 
ran out of memory for degree $n > 7$ (\cite{FongDarve}, Fig.~8),
although the precise memory consumption value was not given.
  
\begin{table}[htb]
\centering
\begin{tabular}{|c|c|c|c|c|c|}
\hline
system size $N$ & 10K & 80K & 640K & 5.12M & 40.96M \\
\hline
direct sum memory (MB) & 0.88 & 7.04 & 56.36 & 450.60 & 3604.80\\
\hline
treecode memory (MB) & & & & & \\
$n=1$ & 0.89 & ~\,7.09 & 56.74 & 453.69 & 3629.24 \\
$n=3$ & 0.91 & ~\,7.27 & 58.17 & 465.09 & 3720.46 \\
$n=5$ & 0.98 & ~\,7.74 & 61.96 & 495.37 & 3962.69 \\
$n=7$ & 1.16 & ~\.8.66 & 69.28 & 553.96 & 4431.40 \\
$n=9$ & 1.44 & 10.16   & 81.32 & 650.30 & 5202.10 \\
\hline
\end{tabular}
\caption{Example 1,
regularized Stokeslets  in a cube~\eqref{Velocity_1},
MRS parameter $\epsilon = 0.02$,
direct sum and KITC memory consumption $\rm{(MB)}$,
system size $N$,
treecode MAC parameter $\theta = 0.7$, degree $n=1,3,5,7,9$, maximum leaf size $N_0 = 2000$.}
\label{tab:TreecodePerformance_Memory}
\end{table}
   
\subsubsection{Effect of MRS parameter}

Table~\ref{table:CaseI_DiffEps} 
shows the error $E$ for values of the MRS parameter 
in the range $0.005 \le \epsilon \le 0.08$,
with microorganism length $l = \epsilon$
and
system size $N$ = 640K~\cite{RostamiOlson}.
The KITC used MAC $\theta = 0.7$
and
degree $n=7$,
and
the error is on the order of 1e-5.
As $\epsilon$ increases,
the error first increases and then decreases;
some error variation with $\epsilon$ is expected,
but the reason for this particular variation is not clear.
It can be noted that error variation with $\epsilon$ also occurred in 
KIFMM simulations~(\cite{RostamiOlson}, Table 3).

\begin{table}[htb]
\centering 
\begin{tabular}{|c|c|c|c|c|c|}
\hline
MRS parameter, $\epsilon$ & 0.005 & 0.01 & 0.02 & 0.04 & 0.08 \\
\hline
error, $E$ & 7.68e-6 & 2.08e-5 & 3.17e-5 & 2.78e-5 & 1.93e-5 \\
\hline
\end{tabular}
\caption{
Example 1,
regularized Stokelets  in a cube~\eqref{Velocity_1},
system size $N\,{=}\,$640K,
MRS parameter $\epsilon$,
microorganism length $l = \epsilon$,
error $E$,
treecode parameters $\theta = 0.7, n = 7$, maximum leaf size $N_0 = 2000$.}
\label{table:CaseI_DiffEps}
\end{table}

 
\subsubsection{Parallel simulations}

We parallelized the direct sum and KITC using OpenMP with up to 24 cores on a single node.
In both methods, the loop over target particles was parallelized,
and
the computation of the modified weights in the KITC was also parallelized.

We consider a system of $N\,{=}\,$640K regularized Stokeslets;
the treecode parameters are $\theta = 0.7, n=7, N_0 = 2000$,
yielding $E\,{=}\,$3.17e-5.
Table~\ref{table:case_1_parallel} shows the 
CPU time for direct sum and KITC $(d,t)$, 
parallel speedup and efficiency,
and 
ratio of direct sum and treecode CPU time.
Parallelizing the direct sum reduces the CPU time 
from 2249\,s with 1~core to 108\,s with 24~cores,
yielding parallel efficiency $86\%$.
Parallelizing the KITC reduces the CPU time 
from 149\,s with 1~core to 10\,s with 24~cores,
yielding parallel efficiency $61\%$.
The KITC is 15 times faster than direct summation with 1~core
and 
10~times faster with 24~cores.

For comparison, 
the KIFMM results in~\cite{RostamiOlson}
employed the Matlab Parallel Computing Toolbox with up to 8~cores.
For a system of size $N\,{=}\,$640K using a $4^3$ grid to solve for the equivalent densities,
the error was $E\,{=}\,$7.33e-4,
and
the CPU time was reduced from 537\,s with 1~core to 100\,s with 8~cores,
yielding parallel efficiency $67$\% (\cite{RostamiOlson}, Tables 4 and 5).
In comparing the KITC results and KIFMM results in~\cite{RostamiOlson},
the point is not to claim that one method is better than the other;
that would require more extensive testing beyond the scope of this work,
but we believe the present results indicate that the KITC is a competitive option for 
fast summation of MRS kernels.

\begin{table}[htb] 
\centering 
\begin{tabular}{|c|c c c|c c c|c|}
\hline
$nc$ & 
$d$ time (s) & $d_1/d_{nc}$ & $d$ PE (\%) & 
$t$ time (s) & $t_1/t_{nc}$ & $t$ PE (\%) & $d/t$ \\
\hline
1   & 2249.71 & 1.00       & 100.00~ & 149.75\, & 1.00      & 100.00\,~ & 15.02 \\
2   & 1123.72 & 2.00       & 100.00~ & ~75.01   & 2.00      & 99.82       & 14.98 \\
4   & ~569.58 & 3.95       & 98.74     & ~38.83   & 3.86      & 96.41       & 14.67 \\
8   & ~311.50 & 7.22       & 90.28     & ~21.44   & 6.98       & 87.31      & 14.53 \\
12 & ~215.49 & 10.43\,~ & 87.00     & ~16.34   & 9.16       & 76.37      & 13.19 \\
24 & ~108.22 & 20.79\,~ & 86.62     & ~10.19   & 14.70\,~ & 61.23      & 10.62 \\
\hline
\end{tabular}
\caption{
Example 1,
regularized Stokeslets  in a cube~\eqref{Velocity_1},
parallel performance,
system size $N\,{=}\,$640K,
MRS parameter $\epsilon = 0.02$,
KITC parameters $\theta = 0.7, n = 7, N_0=2000$, error $E$ = 3.17e-5, 
number of cores ($nc$), 
CPU time (s) is shown for direct sum and KITC $(d,t)$,
parallel speedup $= d_1/d_{nc}, t_{1}/t_{nc}$,
parallel efficiency (PE \%) = $(d_1/d_{nc})/nc, (t_1/t_{nc})/nc$,
parallel KITC speedup $(d/t)$.
} 
\label{table:case_1_parallel}
\end{table}



\subsection{Example 2}

\begin{figure}[htb]
\centering
\includegraphics[width=0.75\textwidth]{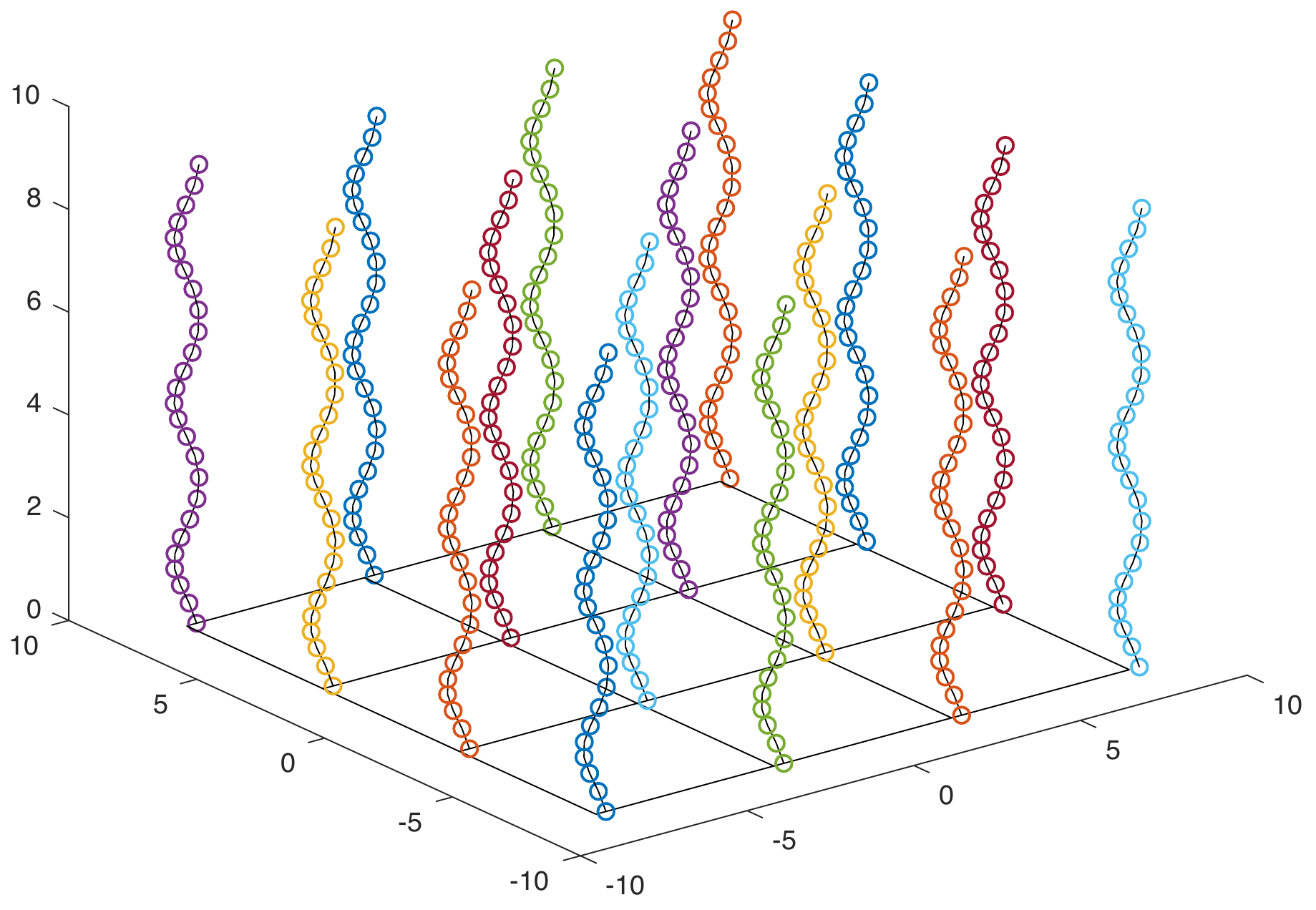}
\caption{Example 2,
regularized Stokeslets and rotlets on an array of helical rods 
representing cilia or free-swimming flagella~\cite{OlsonLimCortez,RostamiOlson}.}
\label{fig:rod}
\end{figure}

The second example models an array of rods representing 
cilia or free-swimming flagella~\cite{OlsonLimCortez,RostamiOlson}.
The number of rods is $N_r$
and
each rod is a helical curve with $M$ segments,
so the total number of particles is $N = N_r(M + 1)$.  
Each rod is parametrized by the $z$-coordinate
and
has the form $(x_0 + 0.3\cos 2z, y_0 + 0.3\sin 2z, z)$ for $0 \le z \le 9$,
where the base point $(x_0,y_0)$ lies at a regular Cartesian grid point in the $xy$-plane
and
the rods extend in the $z$-direction.
Figure~\ref{fig:rod} shows an example with $N_r = 16$ rods
and
$M = 25$ segments on each rod.

In this example each particle is a superposition of a regularized Stokeslet and rotlet,
and
the KITC is applied to compute the linear velocity~\eqref{Velocity_2} 
and 
angular velocity~\eqref{AngularVelocity} for a given rod configuration.
Starting with $N_r = 15^2$ rods with base points in the domain
$[-8,8] \times [-8,8]$ as in~\cite{RostamiOlson},
we increase the number of rods 
($N_r = 15^2, 20^2, 30^2, 40^2, 60^2, 80^2$) while expanding the 
horizontal dimensions of the domain to maintain constant rod density.
Each rod has $M = 150$ segments,
the MRS parameter is $\epsilon = 5L/M = 0.3$,
and
each component of the force ${\bf f}_j$ and torque ${\bf n}_j$ is a random number in $[-1,1]$.

Unlike Example 1 where the particles lie in a cube,
in Example 2 the particles lie in a rectangular slab
(the $z$-direction is shorter than the $xy$-directions).
In particular,
the slab dimensions vary from
$16^2 \times 9$ for the smallest system $(N_r = 15^2, N = 33975)$ to
approximately $85^2 \times 9$ for the largest system $(N_r = 80^2, N=966400)$.
The tree construction scheme described above yields clusters that are
well adapted to the slab geometry of this example.

Also note that the regularized Stokeslet/rotlet kernel~\eqref{MRS_2}
in Example~2 
has more terms 
and
is more expensive to evaluate than the regularized 
Stokeslet kernel~\eqref{Velocity_1} in Example~1.
Hence to better balance the cost of 
direct particle-particle interactions~\eqref{pc_interaction_3D}
and
particle-cluster approximations~\eqref{pc_approximation_3D},
the maximum leaf size was reduced to $N_0 = 1000$ in Example~2.
  
\subsubsection{Error and CPU time versus system size}
 
\begin{figure}[htb]
\centering
\includegraphics[width=0.495\textwidth]{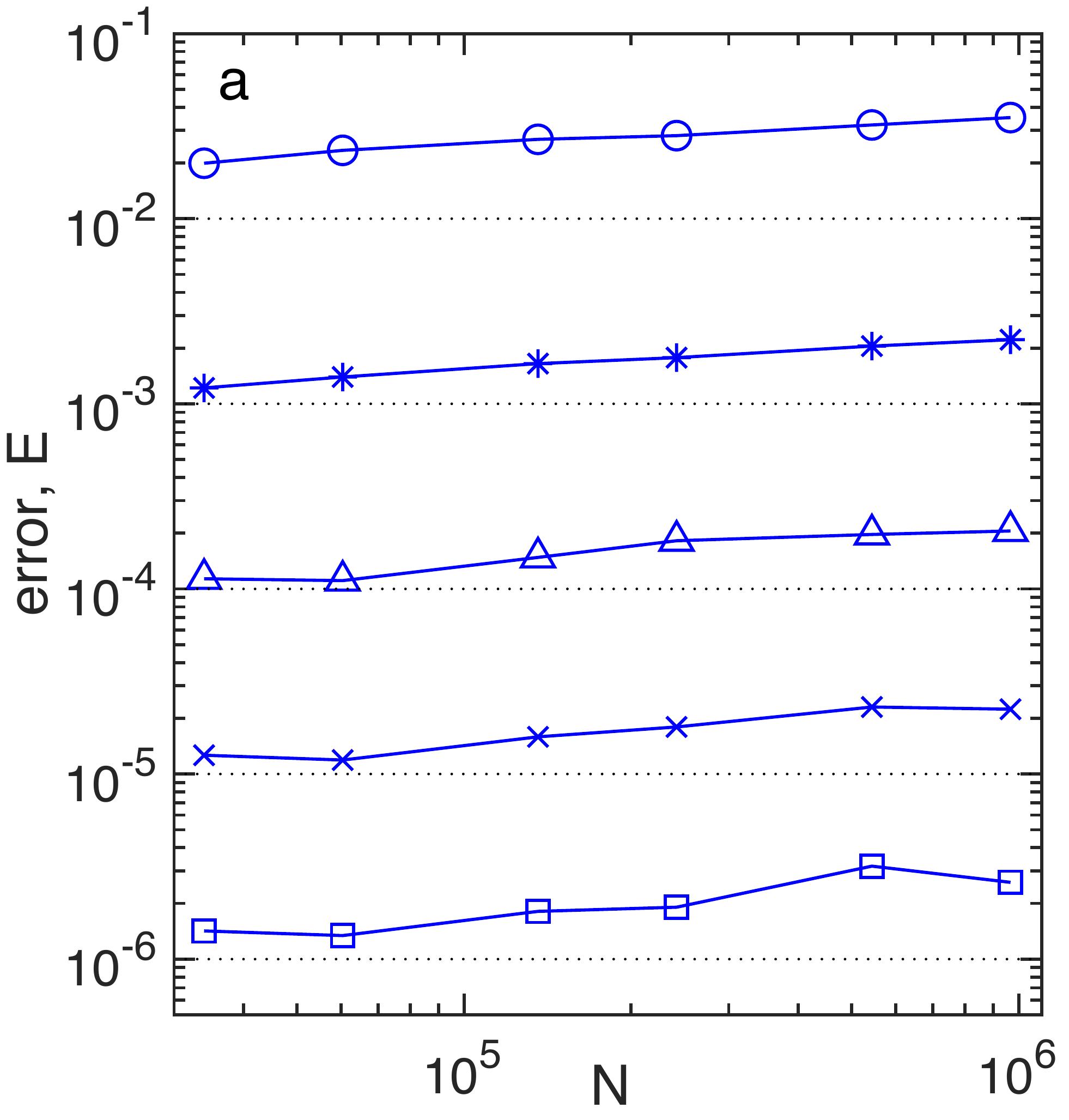}
\includegraphics[width=0.495\textwidth]{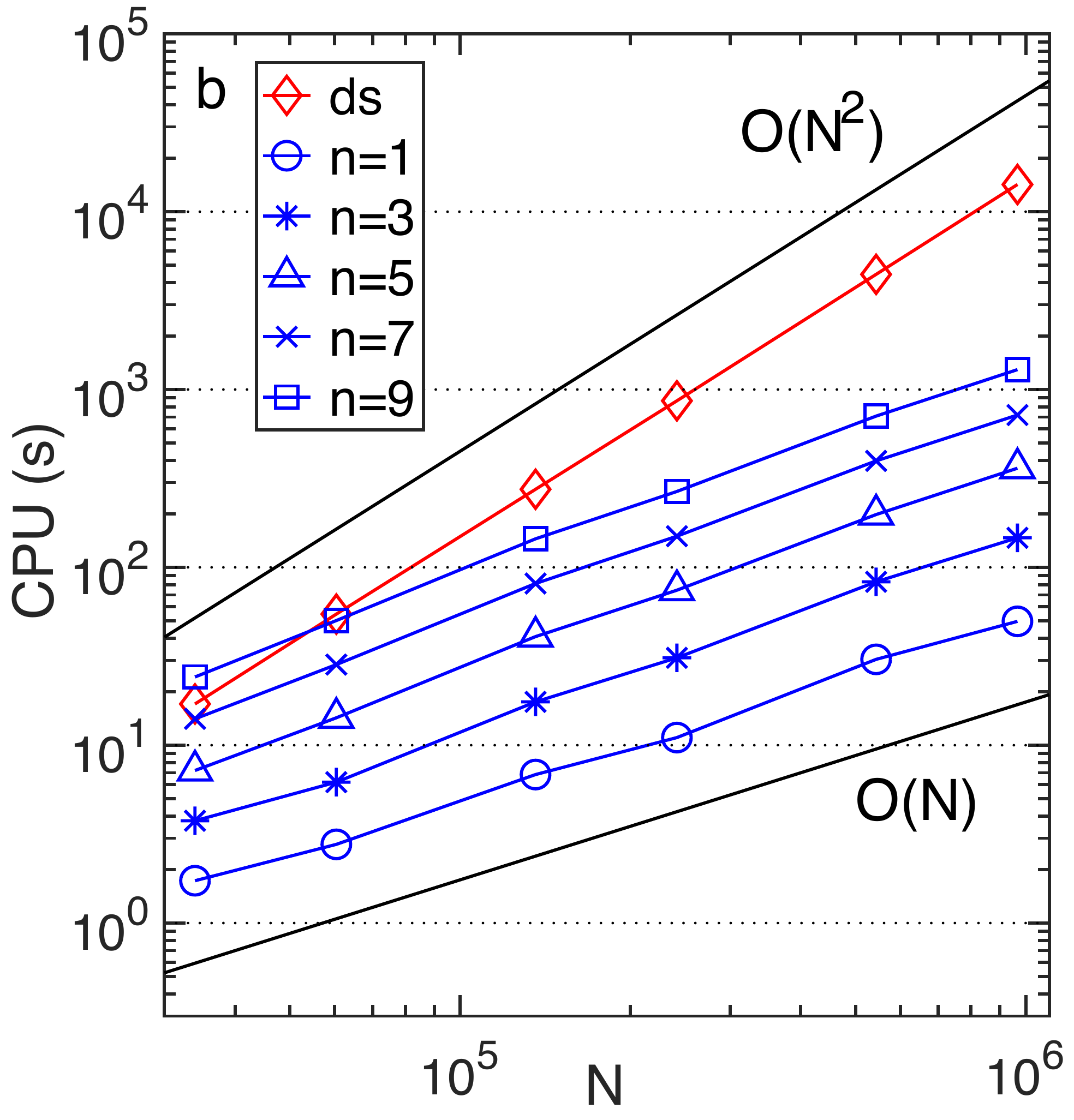}
\vskip -5pt
\caption{Example 2, 
regularized Stokeslets and rotlets on an array of 
helical rods~\eqref{Velocity_2}-\eqref{AngularVelocity}
following~\cite{RostamiOlson},
(a) KITC error $E$,
(b) CPU time~$\rm{(s)}$ for direct sum (red,$\Diamond$)
and KITC (blue, other symbols),
number of rods $N_r=15^2, 20^2, 30^2, 40^2, 60^2, 80^2$,
number of segments per rod $M=150$,
system size $N = N_r(M+1)$,
MRS parameter $\epsilon = 0.3$,
KITC parameters $\theta = 0.7, n = 1,3,5,7,9, N_0 = 1000$.}
\label{fig:CaseII_CPUErr_N}
\end{figure}

Figure~\ref{fig:CaseII_CPUErr_N} presents
the treecode error (a)
and
the direct sum and treecode CPU time (b)
versus system size $N$.
The treecode MAC parameter is $\theta = 0.7$ 
and 
the degree is $n = 1,3,5,7,9$. 
The error increases slightly with system size $N$,
but for a given degree $n$,
the error amplitude is comparable to the results in Example 1 for this range of system size.
As before the direct sum CPU time is parallel to the slope~2 line,
while the KITC time is slightly steeper than the slope~1 line
and is consistent with $O(N\log N)$ scaling.

\subsubsection{Parallel simulations}

Table~\ref{case_2_parallel} presents results 
for $N_r = 80^2$ rods, system size $N = 966400$,
and
with treecode parameters $\theta = 0.7, n=7, N_0 = 1000$, yielding error $E$ = 2.24e-5.
Parallelizing the direct sum reduces the CPU time from
15651\,s with 1~core to 742\,s with 24~cores,
yielding parallel efficiency $87\%$,
while parallelizing the KITC reduces the CPU time from 
720\,s with 1~core to 43\,s with 24~cores, 
yielding parallel efficiency $68\%$. 
The KITC is 21 times faster than direct summation with 1~core 
and 17 times faster with 24~cores.

\begin{table}[htb]
\centering 
\begin{tabular}{|c|ccc|ccc|c|}
\hline
$nc$ & 
$d$ time (s) & $d_1/d_{np}$ & $d$ PE (\%) & 
$t$ time (s) & $t_1/t_{np}$ & $t$ PE (\%) & $d/t$ \\
\hline
1   & 15651.61 & ~~1.00 & 100.00~ & 720.87 & ~~1.00 & 100.0~ & 21.71 \\
2   & ~7907.10 & ~~1.98 & 98.97     & 371.91 & ~~1.94 & 96.91   & 21.26 \\
4   & ~3962.36 & ~~3.95 & 98.75     & 191.89 & ~~3.76 & 93.91   & 20.65 \\
8   & ~2120.36 & ~~7.38 & 92.27     & 110.88 & ~~6.50 & 81.27   & 19.12 \\
12 & ~1424.27 & 10.99   & 91.58     & ~80.35 & ~~8.97 & 74.76   & 17.73 \\
24 & ~~742.85 & 21.07   & 87.79     & ~43.63 & 16.52   & 68.84   & 17.03 \\
\hline
\end{tabular} 
\caption{
Example 2,
regularized Stokeslets and rotlets on an array of helical rods~\eqref{Velocity_2}-\eqref{AngularVelocity}
following~\cite{RostamiOlson},
parallel performance,
$N_r = 80^2$ rods, system size $N = 966400$,  
MRS parameter $\epsilon = 0.3$,
KITC parameters $\theta = 0.7, n = 7$, $N_0 = 1000$, error $E$ = 2.24e-5,
number of cores ($nc$), 
CPU time (s) is shown for direct sum and KITC $(d,t)$, 
parallel speedup = $d_1/d_{nc}, t_1/t_{nc}$, 
parallel efficiency (PE \%) = $(d_1/d_{nc})/nc, (t_1/t_{nc})/nc$,
parallel KITC speedup $(d/t)$.
}
\label{case_2_parallel}
\end{table}

\section{Summary}
\label{sec:summary}

We presented a kernel-independent treecode (KITC) for 
fast summation of particle interactions in 3D.
The method employs barycentric Lagrange interpolation at Chebyshev points 
to approximate well-separated particle-cluster interactions.
The KITC requires only kernel evaluations,
is suitable for non-oscillatory kernels,
and
it utilizes on a scale-invariance property of barycentric Lagrange interpolation.
Numerical results were presented for the non-homogeneous kernels arising in the
Method of Regularized Stokeslets (MRS)~\cite{RostamiOlson}.
For a given level of accuracy, 
the treecode CPU time scales like $O(N \log N )$, 
where $N$ is the number of particles, 
and 
a substantial speedup over direct summation is achieved for large systems. 
The KITC is a relatively simple algorithm with low memory consumption,
and this enables a straightforward parallelization;
here we employed OpenMP with up to 24 cores on a single node;
alternative approaches including distributed memory parallelization for larger systems
will be considered in the future~\cite{feng-barua-li-li-14,marzouk-ghoniem-05,warren-salmon-93}.

Other kernels for which the KITC may be suitable are the 
regularized Green's functions used in computing nearly singular 
integrals~(e.g.~\cite{beale-lai-01,tlupova-beale-13,tlupova-beale-18}),
the RPY (Rotne-Prager-Yamakawa) tensor for hydrodynamic interactions~\cite{Liang},
and
the generalized Born potential for implicit solvent modeling~\cite{xu-cheng-yang-11}.
We expect that the KITC can take advantage of several techniques
employed in advanced implementations of the KIFMM and bbFMM such as
blocking operations, utilizing BLAS routines,
AVX and SSE vectorization, and GPU and Phi coprocessing
(e.g.~\cite{takahashi-cecka-darve-12,agullo-darve-takahashi-14,lashuk-biros-12,malhotra-biros-16}).
It should be noted that an FMM
could be implemented using 
barycentric Lagrange interpolation for both the target and source variables,
and
this is also an interesting direction for future study.
Finally we mention that the extension to barycentric Hermite interpolation
has been carried out for scalar electrostatic 
kernels~\cite{krasny-wang-19,vaughn-wilson-wang-krasny-19}.

\section*{Acknowledgments}
We are grateful for a computer allocation on the
University of Wisconsin-Milwaukee Mortimer Faculty Research Cluster.
We thank Nathan Vaughn and Leighton Wilson for help with programming.
Partial support is acknowledged from start-up funds provided by the 
University of Wisconsin-Milwaukee, 
NSF grant DMS-1819094, 
the Michigan Institute for Computational Discovery and Engineering (MICDE),
and
the Mcubed program at the University of Michigan.



\begin{thebibliography}{9}

\bibitem{agullo-darve-takahashi-14}
E. Agullo, B. Bramas, O. Coulaud, E. Darve, M. Messner and T. Takahashi,
Task-based FMM for multicore architectures,
SIAM J. Sci. Comput., 36 (2014), C66-C93.

\bibitem{Klinteberg2017}
L. af Klinteberg, D. S. Shamshirgar and A.-K. Tornberg,
Fast Ewald summation for free-space Stokes potentials,
Res. Math. Sci., 4 (2017), Article 1. 


\bibitem{anderson-92}
C. R. Anderson,
An implementation of the Fast Multipole Method without multipoles,
SIAM J. Sci. Stat. Comput., 13 (1992), 923-947.

\bibitem{BarnesHut}
J. E. Barnes and P. Hut,
A hierarchical $O(N\log N)$ force-calculation algorithm,
Nature, 324 (1986), 446-449.

\bibitem{beale-lai-01}
J. T. Beale and M.-C. Lai,
A method for computing nearly singular integrals, 
SIAM J. Numer. Anal., 38 (2001), 1902-1925.

\bibitem{BerrutTrefethen}
J.-P. Berrut and L. N. Trefethen,
Barycentric Lagrange interpolation,
SIAM Rev., 46 (2004), 501-517.

\bibitem{borm-grasedyck-hackbusch-03}
S. B\"orm, L. Grasedyck and W. Hackbusch,
Introduction to hierarchical matrices with applications,
Eng. Anal. Bound. Elem., 27 (2003), 405-422.
 
\bibitem{BouzarthMinion}
{\sc E. L. Bouzarth and M. L. Minion},
{\it Modeling slender bodies with the method of regularized Stokeslets} 
J. Comput. Phys., 230 (2011), pp.~3929--3947.


\bibitem{brandt-lubrecht-90}
A. Brandt and A. A. Lubrecht,
Multilevel matrix multiplication and fast solution of integral equations,
J. Comput. Phys., 90 (1990), 348-370.


\bibitem{cheng-greengard-rokhlin-99}
H. Cheng, L. Greengard and V. Rokhlin,
A fast adaptive multipole algorithm in three dimensions,
J. Comput. Phys., 155 (1999), 468-498.

\bibitem{Cortez}
R. Cortez,
The method of regularized Stokeslets,
SIAM J. Sci. Comput., 23 (2001), 1204-1225.

\bibitem{CortezFauciMedovikov}
R. Cortez, L. Fauci  and A. Medovikov,
The method of regularized Stokeslets in three dimensions: 
Analysis, validation, and application to helical swimming,
Phys. Fluids, 17 (2005), Article 031504.


\bibitem{Draghicescu}
C. I. Draghicescu and M. Draghicescu,
A fast algorithm for vortex blob interactions,
J. Comput. Phys., 116 (1995), 69-78.

\bibitem{driscoll-hale-trefethen-14}
T. A. Driscoll, N. Hale and L. N. Trefethen,
Chebfun Guide, Pafnuty Publications, Oxford, 2014.
www.chebfun.org

\bibitem{duan-k-01}
Z.-H. Duan and R. Krasny,
An adaptive treecode for computing nonbonded potential energy in classical molecular systems, 
J. Comput. Chem., 22 (2001), 184-195.

\bibitem{Essmann}
U. Essmann, L. Perera, M. Berkowitz, T. Darden, H. Lee and L. Pedersen,
A smooth particle mesh Ewald method,
J. Chem. Phys., 103 (1995), 8577-8593. 

\bibitem{feng-barua-li-li-14}
H. Feng, A. Barua, S. Li and X. Li,
A parallel adaptive treecode algorithm for evolution of elastically stressed solids,
Commun. Comput. Phys., 15 (2014), 365-387.

\bibitem{Flores}
H. Flores, E. Lobaton, S. M\'endez-Diez, S. Tlupova, and R. Cortez,
A study of bacterial flagellar bundling, 
Bull. Math. Biol., 67 (2005), 137-168.

\bibitem{FongDarve}
W. Fong and E. Darve,
The black-box fast multipole method,
J. Comput. Phys., 228 (2009), 8712-8725.

\bibitem{GengKrasny}
W.-H. Geng and R. Krasny, 
A treecode-accelerated boundary integral Poisson-Boltzmann solver for solvated biomolecules,
J. Comput. Phys., 247 (2013), 62-78.

\bibitem{giebermann-01}
K. Giebermann,
Multilevel approximation of boundary integral operators,
Computing, 67 (2001), 183-207.

\bibitem{GimbutasRokhlin}
Z. Gimbutas and V. Rokhlin,
A generalized Fast Multipole Method for nonoscillatory kernels,
SIAM J. Sci. Comput., 24 (202), 796-817.



\bibitem{GreengardHuang}
L. F. Greengard and J. Huang,
A new version of the Fast Multipole Method for screened Coulomb interactions in three dimensions,
J. Comput. Phys., 180 (2002), 642-658.

\bibitem{GreengardRokhlin} 
L. Greengard and V. Rokhlin,
A fast algorithm for particle simulations,
J. Comput. Phys., 73 (1987), 325-348.

\bibitem{Greengard}
L. Greengard,
The Rapid Evaluation of Potential Fields in Particle Systems,
MIT Press, Cambridge, MA 1988.

\bibitem{hackbusch-nowak-89}
W. Hackbusch and Z. P. Nowak,
On the fast matrix multiplication in the boundary element method by panel clustering,
Numer. Math., 54 (1989), 463-491.

\bibitem{hardy-skeel-16}
D. J. Hardy, M. A. Wolff, J. Xia, K. Schulten and R. D. Skeel,
Multilevel summation with B-spline interpolation for pairwise interactions
in molecular dynamics simulations,
J. Chem. Phys., 144 (2016), Article 114112.

\bibitem{Higham}
N. J. Higham,
The numerical stability of barycentric Lagrange interpolation,
IMA J. Numer. Anal., 24 (2004), 547-556.

\bibitem{HockneyEastwood}
R. W. Hockney and J. W. Eastwood,
Computer Simulation Using Particles,
Taylor \& Francis, Bristol, 1988.


\bibitem{KrasnyWang}
R. Krasny and L. Wang,
Fast evaluation of multiquadric RBF sums by a Cartesian treecode,
SIAM J. Sci. Comput., 33 (2011), 2341-2355. 

\bibitem{krasny-wang-19}
R. Krasny and L. Wang,
A treecode based on barycentric Hermite interpolation for electrostatic particle interactions,
submitted to proceedings of NSF-CBMS Conference:~Molecular Bioscience and Biophysics,
University of Alabama, Tuscaloosa, May 13-17, 2019

\bibitem{kress-14}
R. Kress,
Linear Integral Equations,
Springer, New York, 2014, third edition.


\bibitem{lashuk-biros-12}
I. Lashuk, A. Chandramowlishwaran, H. Langston, T.-A. Nguyen, R. Sampath, 
A. Shringarpure, R. Vuduc, L. Ying, D. Zorin and G. Biros,
A massively parallel adaptive fast-multipole method on heterogeneous architectures,
Commun. ACM, 55 (2012), 101-109.

\bibitem{letourneau-cecka-darve-14}
P.-D. L\'etourneau, C. Cecka and E. Darve,
Cauchy Fast Multipole Method for general analytic kernels,
SIAM J. Sci. Comput., 36 (2014), A396-A426.

\bibitem{LiJohnstonKrasny}
P. Li, H. Johnston and R. Krasny, 
A Cartesian treecode for screened Coulomb interactions,
J. Comput. Phys., 228 (2009), 3858-3868.

\bibitem{Liang}
Z. Liang, Z. Gimbutas, L. Greengard, J. Huang and S. Jiang,
A fast multipole method for the Rotne-Prager-Yamakawa tensor and its applications,
J. Comput. Phys., 234 (2013), 133-139.

\bibitem{lindsay-krasny-01}
K. Lindsay and R. Krasny,
A particle method and adaptive treecode for vortex sheet motion in three-dimensional flow,
J. Comput. Phys., 172 (2001), 879-907.


\bibitem{makino-99}
J. Makino,
Yet another fast multipole method without multipoles - Pseudoparticle multipole method,
J. Comput. Phys., 151 (1999), 910-920.


\bibitem{malhotra-biros-16}
D. Malhotra and G. Biros,
Algorithm 967: A distributed-memory fast multipole method for volume potentials,
ACM Trans. Math. Softw., 43 (2016), Article 17.

\bibitem{march-xiao-biros-15}
W. B.~March, B. Xiao and G. Biros,
ASKIT:~Approximate skeletonization kernel-independent treecode in high dimensions, 
SIAM J. Sci. Comput., 37 (2015), A1089-A1110.

\bibitem{martinsson-rokhlin-07}
P. G. Martinsson and V. Rokhlin,
An accelerated kernel-independent fast multipole method in one dimension,
SIAM J. Sci. Comput., 29 (2007), 1160-1178.

\bibitem{mascarenhas-14}
W. F. Mascarenhas,
The stability of barycentric interpolation at the Chebyshev points of the second kind,
Numer. Math., 128 (2014), 265-300.

\bibitem{marzouk-ghoniem-05}
Y. M. Marzouk and A. F. Ghoniem,
$K$-means clustering for optimal partitioning and dynamic load balancing of 
parallel hierarchical $N$-body simulations,
J. Comput. Phys., 207 (2005), 493-528.



\bibitem{OlsonLimCortez}
S. D. Olson, S. Lim and R. Cortez,
Modeling the dynamics of an elastic rod with intrinsic curvature and twisting using a 
regularized Stokes formulation,
J. Comput. Phys., 238 (2013), 169-187.

\bibitem{Pozrikidis}
C. Pozrikidis,
Boundary Integral and Singularity Methods for Linearized Viscous Flow,
Cambridge University Press, Cambridge, UK, 1992.

\bibitem{rack-reimer-82}
H.-J. Rack and M. Reimer,
The numerical stability of evaluation schemes for polynomials based on the Lagrange interpolation form,
BIT, 22 (1982), 101-107.

\bibitem{RostamiOlson}
M. W. Rostami and S. D. Olson,
Kernel-independent fast multipole method within the framework of regularized Stokeslets,
J. Fluid Struct., 67 (2016), 60-84.




\bibitem{Salzer}
H. E. Salzer,
Lagrangian interpolation at the Chebyshev points $x_{n,\nu} = \cos(\nu\pi/n), \nu = 0(1)n$; 
some unnoted advantages,
Comput. J., 15 (1972), 156-159.


\bibitem{Smith}
D. J. Smith,
A boundary element regularized Stokeslet method applied to cilia- and flagella-driven flow,
Proc. R. Soc. A, 465 (2009), 3605-3626.

\bibitem{takahashi-cecka-darve-12}
T. Takahashi, C. Cecka and E. Darve,
Optimization of the parallel black-box fast multipole method on CUDA,
in Proceedings of Conference on Innovative Parallel Computing (InPar), 2012,
San Jose, CA, USA. DOI:~10.1109/InPar.2012.6339607.


\bibitem{tausch-03}
J. Tausch,
The Fast Multipole Method for arbitrary Green's functions,
Contemp. Math., 329 (2003), 307-314.

\bibitem{TlupovaCortez}
S. Tlupova and R. Cortez,
Boundary integral solutions of coupled Stokes and Darcy flows,
J. Comput. Phys., 228 (2009), 158-179.

\bibitem{tlupova-beale-13}
S. Tlupova and J. T. Beale,
Nearly singular integrals in 3D Stokes flow,
Commun. Comput. Phys., 14 (2013), 1207-1227.

\bibitem{tlupova-beale-18}
S. Tlupova and J. T. Beale,
Regularized single and double layer integrals in 3D Stokes flow,
J. Comput. Phys., 386 (2019), 568-584.

\bibitem{trefethen-13}
L. N. Trefethen,
Approximation Theory and Approximation Practice,
SIAM, Philadelphia, 2013.


\bibitem{vaughn-wilson-wang-krasny-19}
N. Vaughn, L. Wilson, L. Wang and R. Krasny, 
GPU-accelerated barycentric treecodes,
submitted to proceedings of SIAM Conference on Parallel Processing 2020

\bibitem{WangTlupovaKrasny}
L. Wang, S. Tlupova and R. Krasny,
A treecode algorithm for 3D Stokeslets and stresslets,
Adv. Appl. Math. Mech., 11 (2019), 737-756.

\bibitem{warren-salmon-93}
M. S. Warren and J. K. Salmon,
A parallel hashed oct-tree N-body algorithm,
in Proceedings of the 1993 ACM/IEEE Conference on Supercomputing, (1993), 12-21.




\bibitem{xu-cheng-yang-11}
Z. Xu, X. Cheng and H. Yang,
Treecode-based generalized Born method,
J. Chem. Phys., 134 (2011), Article 064107.

\bibitem{YingBirosZorin}
L. Ying, G. Biros and D. Zorin,
A kernel-independent adaptive fast multipole algorithm in two and three dimensions,
J. Comput. Phys., 196 (2004), 591-626. 

\bibitem{Ying}
L. Ying,
A kernel independent fast multipole algorithm for radial basis functions,
J. Comput. Phys., 213 (2006), 451-457.

\bibitem{zhang-sun-11}
B. Zhang, J. Huang, N. P. Pitsianis and X. Sun,
A Fourier-series-based kernel-independent fast multipole method,
J. Comput. Phys., 230 (2011), 5807-5821.

\end{thebibliography}
\end{document}